\documentclass[twoside,12pt]{article}
\usepackage[]{amsmath,amsthm,amssymb}
\usepackage[latin1]{inputenc}

\begin{document}
\title{{\Large\bf  Certain identities, connection  and explicit formulas for the Bernoulli, Euler  numbers\\  and Riemann  zeta -values}}

\author{Semyon   YAKUBOVICH}
\maketitle

\markboth{\rm \centerline{ Semyon YAKUBOVICH}}{}
\markright{\rm \centerline{Bernoulli, Euler  numbers and Riemann  zeta - values}}

\begin{abstract} {Various new identities,  recurrence relations,   integral representations, connection  and explicit formulas
are established for the Bernoulli, Euler  numbers and the values of
Riemann's   zeta function $\zeta (s)$. To do this, we explore
properties of some Sheffer's  sequences of polynomials related to
the Kontorovich-Lebedev transform. }

\end{abstract}
\vspace{4mm}



\vspace{4mm}

\section {Introduction and preliminary results}

In 2009 the author \cite{yakpol} introduced a family  of polynomials of degree $n$, which belongs to Sheffer's sequences (cf. \cite{bro}) and related to the Kontorovich-Lebedev transform.  Precisely, it has the form
$$p_n(x)= (-1)^n e^x {\cal A} ^n e^{-x}, \quad n \in \mathbb{N}_0,\eqno(1.1)$$
where
$${\cal A} \equiv x^2- x{d\over dx} x {d\over dx},\eqno(1.2)$$
is the second order differential operator having as an eigenfunction the modified Bessel function $ K_{i\tau}(x),\  \tau \in \mathbb{R},\ $ ($ i$  is the imaginary unit ), i.e.
$${\cal A} \ K_{i\tau}(x)= \tau^2 K_{i\tau}(x).\eqno(1.3)$$
The modified Bessel function $ K_{i\tau}(x)$ is, in turn,  the kernel of the Kontorovich-Lebedev transform (see in
\cite{sned}, \cite{leb},  \cite{yakind})
$$(Gf) (\tau)=\int_0^\infty K_{i\tau}(x)f(x) {dx\over x}, \ \tau \in \mathbb{R}_+.\eqno(1.4)$$
As it is known,  operator (1.4) extends to a bounded invertible map  $G: L_2\left({\mathbb {R}}_+;  x^{-1}dx  \right)
\to  L_2({\mathbb {R}}_+;\tau\sinh\pi\tau d\tau)$ and this map is isometric, i.e.
$$\int_0^\infty \tau\sinh\pi\tau |(Gf)(\tau)|^2 d\tau= {\pi^2\over 2} \int_0^\infty |f(x)|^2 {dx\over x} .\eqno(1.5)$$
Reciprocally,  the inversion formula holds
$$f(x)={2\over  \pi^2}\int_{0}^\infty \tau\sinh\pi\tau K_{i\tau}(x)(Gf)(\tau) d\tau,\  x >0.\eqno(1.6)$$
The modified Bessel function $K_{i\tau}(x)$ has  the asymptotic behavior \cite{leb}
$$ K_\nu(z) = \left( \frac{\pi}{2z} \right)^{1/2} e^{-z} [1+
O(1/z)], \qquad z \to \infty,\eqno(1.7)$$
and near the origin
$$ K_\nu(z) = O\left ( z^{-|{\rm Re}\nu|}\right), \ z \to 0,\eqno(1.8)$$
$$K_0(z) = -\log z + O(1), \ z \to 0. \eqno(1.9)$$
Moreover it can be defined by the  following integral
representations
$$K_{\nu}(x)=\int_0^\infty e^{-x\cosh u}\cosh \nu u du,
x>0, \eqno(1.10)$$
$$K_{\nu}(x)={1\over 2} \left(x\over 2\right)^\nu\int_0^\infty
e^{-t-{x^2\over 4t}}t^{-\nu-1}dt, \ x>0.\eqno(1.11)$$
Returning to the system of polynomials (1.1), we easily observe that all their coefficients $a_{n,k},\  k= 0,1\dots, n$ are integers.  It can be represented by the integral
$$p_n(x)= {2(-1)^n\over \pi} e^{x}\int_{0}^\infty \tau^{2n} K_{i\tau}(x)\ d\tau,\eqno(1.12) $$ and satisfies  the
differential recurrence relation of the form
$$p_{n+1}(x)= x^2p_n^{\prime\prime}(x) + x(1-2x)p_n^\prime(x)- xp_n(x), \ n= 0,1,2,\dots \ .\eqno(1.13)$$
In particular, we derive
$$p_0(x)=1, \ p_1(x)= -x,\ p_2(x)= 3x^2-x, \quad p_3(x)= -15x^3+ 15x^2-x.$$
The leading coefficient $a_{n,n}$ of these polynomials can be calculated by the formula
$$a_{n,n}= (-1)^n (2n-1)!! =  (-1)^n 1\cdot 3\cdot 5\dots \  \cdot (2n-1), \ n \in \mathbb{N}.\eqno(1.14)$$
Moreover, recently we found the explicit formula of coefficients $a_{n,k} $  (see \cite{yakyor})
$$a_{n,k}= {1\over k!}\sum_{r=0}^{k} {(-1)^r\over 2^r} \binom{k}{r}\sum_{j=0}^{k-r} {(-1)^j\over 2^j} \binom{k-r}{j}
(r-j)^{2n} ,\  k= 1,\dots, n\eqno(1.15)$$
and as a consequence of the definition (1.1) $a_{n,k} \in \mathbb{Z}$, i.e. the right- hand side of (1.15) is an integer.
The generating function $\Phi(x,t)$ for this sequence of polynomials is given by the series
$$\Phi(x,t)=  e^{- 2x\sinh^2( t/2)} = \sum_{n=0}^\infty {p_n(x)\over (2n)!}  \ t^{2n},\ |t| < {\pi\over 4}.\eqno(1.16)$$
Letting $x=0$ in the latter equation, we find
$$p_n(0)= 0, \quad n=1,2,\dots.$$
A differentiation with respect to  $x$ in (1.16) yields the equality
$${\partial \Phi \over \partial x } =  (1-\cosh t)\  e^{- 2x\sinh^2( t/2)} = \sum_{n=0}^\infty {p_n^\prime (x)\over (2n)!} \ t^{2n}.\eqno(1.17)$$
Decomposing the left-hand side of (1.17) as a product of series and equating coefficients in front of $t^{2n}$ we come up with  the following recurrence relation
$$p_n^\prime (x)= - \sum_{k=0}^{n-1} {2n \choose 2k} p_k(x), \ n \in
\mathbb{N}.\eqno(1.18)$$
Putting $x=0$ in (1.17) and using  values $p_0(0)= 1, \  p_k(0)=0, \ k \in \mathbb{N}$ we obtain $p_n^\prime (0)= -1, \ n \in  \mathbb{N}$. Analogously, a differentiation with respect to $t$ in  (2.4) leads us to
$${\partial \Phi \over \partial t } =  - x \sinh t \  e^{- 2x\sinh^2( t/2)} = \sum_{n=1}^\infty {p_n (x)\over (2n -1)!}
\ t^{2n -1}.\eqno(1.19)$$
Similarly we derive the relation
$$p_{n+1} (x)= - x \sum_{k=0}^{n} {2n +1 \choose 2k} p_k(x), \ n \in
\mathbb{N}_0.\eqno(1.20)$$
Moreover, differentiating through in (1.19),  we call  (1.17) and using simple
relations for binomial coefficients,  we obtain  the identity
$$ x \sum_{k=0}^{n} {2n +1 \choose 2k} p_k^\prime(x) =  \sum_{k=1}^{n} {2n +1 \choose 2k -1}
p_k(x).\eqno(1.21)$$
Comparing  (1.17), (1.19), we find that $\Phi(x,t)$ satisfies the following first order partial differential equations
$$ {\partial \Phi \over \partial  t } +  x\sinh t\   \Phi = 0,\eqno(1.22)$$
$${\partial \Phi \over \partial x }  +  2   \sinh^2 \left(t\over 2\right)   \Phi =0,\eqno(1.23) $$
$$x {\partial \Phi \over \partial x }  =  \tanh \left(t\over 2\right)   {\partial \Phi \over \partial t }.\eqno(1.24) $$
Further, returning to representation (1.12) and employing the inversion formula (1.6) of the Kontorovich-Lebedev transform, we obtain the equality
$${\tau^{2n-1}\over \sinh\pi\tau} = {(-1)^n\over \pi}\int_0^\infty e^{-x}K_{i\tau}(x)p_n(x){dx\over
x}, \ n \in \mathbb{N}.\eqno(1.25)$$
In particular, it yields
$$\int_0^\infty e^{-x}K_{0}(x)p_n(x){dx\over x} =0, \ n=2,3,\dots \  .$$
Integrating with respect to $\tau$ in (1.25), we call the value of the integral  (2.4.3.1) in \cite{prud},  Vol. I
$$\int_0^\infty {\tau^{\alpha -1} d\tau \over \sinh\pi\tau} =    \frac{2^\alpha - 1}{\pi^\alpha
2^{\alpha-1}}\  \Gamma(\alpha) \zeta(\alpha),\quad {\rm Re} \  \alpha >  1,\eqno(1.26)$$
where $\Gamma(\alpha)$ is Euler's gamma function and $\zeta(\alpha) $ is Riemann's zeta function (cf. \cite{erd}, Vol. 1)
and relation (2.16.48.1) in \cite{prud}, Vol. II to obtain the following representations of zeta-values at  even and odd integers, respectively,
$$\frac{2^{2n}  - 1}{2^{2(n-1)}}\  (-1)^n  (2n-1)! { \zeta(2n)\over \pi^{2n}}  =  \int_0^\infty e^{-2x} p_n(x){dx\over
x}, \ n \in \mathbb{N},\eqno(1.27)$$
$$(-1)^n  (2n)!  \  \left(2^{2n+1} - 1\right)\   {\zeta(2n+1)\over (2 \pi)^{2n} } =  \int_0^\infty \int_0^\infty \tau K_{i\tau} (x) e^{-x} p_n(x){d\tau dx\over  x}, \ n \in \mathbb{N}.\eqno(1.28)$$
Finally in this section we note that recently in \cite{lour} the family  of polynomials (1.1) was generalized on the sequence
$$p_n(x; \alpha)=   (-1)^n e^x x^{-\alpha} {\cal A} ^n e^{-x} x^\alpha , \quad n \in \mathbb{N}_0,\eqno(1.29)$$
involving an arbitrary parameter $\alpha, \  {\rm Re} \  \alpha >  - 1/2. $

\section{Identities for the Bernoulli and Euler numbers}

In this section we will derive a number of recurrence relations,  finite sum, connection and explicit formulas,  series and integral representations, which are related to the Bernoulli and Euler numbers. The Bernoulli numbers $B_n,\  n= 0,1,2,\dots, $ can be defined via the generating function (see in \cite{erd},  Vol. I )
$${x\over e^x- 1} = \sum_{n=0}^\infty B_n {x^n\over n! }, \quad |x| < 2\pi \eqno(2.1)$$
and Bernoulli polynomials $B_n(x)$ by the equality
$${t e^{xt}   \over e^t - 1} = \sum_{n=0}^\infty B_n(x)  {t^n\over n! }, \quad |x| < 2\pi. \eqno(2.1)$$
In particular, we find the values,   $B_0=1,\ B_1= -1/2, B_2= 1/6, B_4= -1/30$ and $B_n=0$ for all odd  $n \ge 3$. Furthermore,   $(-1)^{n-1} B_{2n} > 0$ for all $ n \ge 1$.    We list some important properties of the Bernoulli numbers and polynomials, which will be employed below.  All details and proofs can be found in \cite{erd}, Vol. I.   The basic identity for Bernoulli numbers is
$$ \sum_{k=0}^{n-1} {n \choose k} B_k = 0, \ n \ge 2.\eqno(2.3)$$
Concerning the Bernoulli polynomials,  it has the explicit formula,
$$B_n(x)=   \sum_{k=0}^{n} {n \choose k} B_k x^{n-k}.\eqno(2.4)$$
Hence,  for instance,
$$B_0(x)= 1, \  B_1(x)= x- {1\over 2},\  B_2(x)=  x^2- x+ {1\over 6},$$$$    B_3(x)=  x^3- {3\over 2} x^2+ {1\over 2} x,\
B_4(x) = x^4- 2x^3+ x^2- {1\over 30}.\eqno(2.5)$$
Evidently,  $B_n= B_n(0)$.   The Bernoulli polynomials and their  derivative satisfy the following important relations
$$B_n^\prime (x)=  n B_{n-1} (x),\eqno(2.6)$$
$$B_n(x+1)- B_n(x)= n x^{n-1},\eqno(2.7)$$
$$B_n(1-x)= (-1)^n B_n(x), \eqno(2.8)$$
$$B_n(mx)= m^{n-1} \sum_{k=0}^{m-1} B_n \left(x+  {k\over m}\right).\eqno(2.9)$$
The remarkable Euler formula relates the Bernoulli numbers $B_{2n}$ and Riemann zeta- values $\zeta(2n)$  (see, for instance, in \cite{dwi} )
$$\zeta(2n)=  (-1)^{n-1} \frac{2^{2n-1} B_{2n}}{(2n)!} \pi^{2n}.\eqno(2.10)$$
The Euler numbers $E_n,\  n= 0,1,2,\dots, $ can be defined, in turn,  by the equality  (see in \cite{erd},  Vol. I )
$${1\over \cosh t } = \sum_{n=0}^\infty E_n  {t^n\over n! }, \quad |t| <  {\pi\over 2}. \eqno(2.11)$$
As we see,  $E_{2n+1} = 0$ and, in particular,  $E_0=1,\ E_2=-1, \  E_4=5,\  E_6= -61$.  The basic relation  for Euler  numbers is
$$ \sum_{k=0}^{n} {2n \choose 2k} E_{2k} = 0, \ n \ge 1.\eqno(2.12)$$

In order to obtain new properties of the Bernoulli and Euler numbers, we will employ Sheffer's sequences, which are associated with the Kontorovich-Lebedev transform (1.4).  Indeed,  calling identity (1.27), we immediately derive the integral representation of Bernoulli's numbers $B_{2n}$ in terms of the sequence of polynomials (1.1), namely
$$B_{2n}= \frac{n}{ 1- 2^{2n}} \int_0^\infty e^{-2x} p_n(x){dx\over  x}.\eqno(2.13)$$
For the numbers $B_{4n}$ we have the formula (see in \cite{yakpol})
$$B_{4n} =  \frac{2n}{ 1- 2^{4n}} \int_0^\infty e^{-2x} p^2_n(x){dx\over  x},\eqno(2.14)$$
which is the result of the equality
$$ \int_0^\infty e^{-2x} p_n(x){dx\over  x} =  \int_0^\infty e^{-2x} p^2_n(x){dx\over  x}.\eqno(2.15)$$
But as it is proved in \cite{yakpol},  a more general relation takes place
$$ \int_0^\infty e^{-2x} p_{n+m}(x){dx\over  x} =  \int_0^\infty e^{-2x} p_n(x) p_m(x) {dx\over  x},\eqno(2.16)$$
which holds for any $n, m \in \mathbb{N}_0$ such that at least one is nonzero.   Hence, appealing to (2.13), we derive the identity
$$B_{2(n+m)}= \frac{n+m}{ 1- 2^{2(n+m)}} \int_0^\infty e^{-2x} p_{n-k}(x) p_{m+k}(x) {dx\over  x}\eqno(2.17)$$
being valid for any $k=0,1,\dots, n$.

Another definition of Euler's numbers can be given in terms of  Sheffer's  sequence of polynomials $q_n(x)$ introduced in \cite{yakpol}
$$q_n(x)= e^x\int_x^\infty e^{-t} p_n(t) dt,\quad   n \in \mathbb{N}_0,\eqno(2.18)$$
which, in turn,   is defined via the generating function $F(x,t)= \Phi(x,t) / \cosh t$  (see (1.16))
$${1\over \cosh t} e^{- 2x\sinh^2( t/2)} = \sum_{n=0}^\infty {q_n(x)\over (2n)!}  \ t^{2n},\ |t| < {\pi\over 4}.\eqno(2.19)$$
Hence with the use of (2.11) and the integral representation of Euler numbers \cite{yakpol}
$$E_{2n}=  \int_0^\infty e^{- x} p_n(x) dx\quad  \in \mathbb{Z}\eqno(2.20)$$
 we find
$$E_{2n} = q_n(0), \quad n=0,1,\dots.$$
Moreover, following \cite{yakpol}, the sequence $q_n(x)$ has a relationship with $p_n(x)$.   Indeed,
$$q_n(x)= \sum_{k=0}^n p_n^{(k)}(x),\ n \in \mathbb{N}_0,\eqno(2.21)$$
$$q_n(x)= \sum_{k=0}^n E_{2(n-k)} {2n  \choose 2k }p_{k}(x),\eqno(2.22)$$
where $p_n^{(k)}(x)$ is the $k$ -th derivative of $p_n(x)$ and relation (2.22) can be obtained employing (2.18) and the binomial type identity for the sequence $p_n(x)$ (cf. \cite{bro})
$$p_n(x+y)=  \sum_{k=0}^n {2n \choose 2k}  p_k(x)p_{n-k}(y).\eqno(2.23)$$
For  instance,
$$q_0(x)= 1, \quad  q_1(x)= - (x+1), $$
$$\ q_2(x)= 3x^2 + 5x +5, \quad   q_3(x)= -15x^3-30x^2-61x -61.$$
Differentiating through in (2.21), we find
$$q_n^\prime (x)= \sum_{k=0}^n p_n^{(k+1)}(x)= \sum_{k=1}^{n} p_n^{(k)}(x)= q_n(x)- p_n(x).$$
Thus, $p_n(x)= q_n(x)- q_n^\prime (x)$ and since $p_n(0)=0, n \in
\mathbb{N}$, we get
$$E_{2n}= q_n^\prime (0),\ n \in \mathbb{N}.$$
Moreover,
$$p_n^\prime (x)= q_n^\prime (x)- q_n^{\prime\prime} (x)$$
and since $p_n^\prime (0)= -1, n = 1,2,\dots \ $ (see (1.18)),  it yields
$$q_n^{\prime\prime} (0)= E_{2n}+1,\  n \ge 1.\eqno(2.24)$$
Returning to (2.19), we easily derive analogs of the first order partial differential equations (1.22), (1.23), (1.24) for the generating  function $F(x,t)$,  namely
$$ {\partial F \over \partial  t } +  \left[ x\sinh t\  + \tanh t  \right] F = 0,\eqno(2.25)$$
$${\partial F \over \partial x }  +  2   \sinh^2 \left(t\over 2\right)   F =0,\eqno(2.26) $$
$$x {\partial F \over \partial x }  =  \tanh \left(t\over 2\right)  \left[  {\partial F \over \partial t } + F \tanh t \right].\eqno(2.27) $$
To find the  inverse of relation (2.22), we apply the product of series of  $\cosh t$ and  (2.19), equating its coefficients in front of $t^{2n}$ with the corresponding terms of  series (1.16). This is indeed allowed within the   interval of the absolute convergence $|t| < \pi/4$.  As a result, we deduce
$$p_n(x)= \sum_{k=0}^n  {2n  \choose 2k }q_{k}(x).\eqno(2.28)$$
But $q_n^\prime (x)=  q_n(x)- p_n(x)$.   So, we have
$$q_n^\prime (x) = -  \sum_{k=0}^{n-1}   {2n  \choose 2k }q_{k}(x).\eqno(2.29)$$
Another source of identities for Bernoulli numbers is a formula related to $p^\prime_n(x)$. To derive it, we call the partial differential equation (1.24) and the Taylor series for the hyperbolic tangent
$$\tanh\left({x\over 2}\right)=   2 \sum_{k=0}^\infty \frac{\left(2^{2(k+1)}  -  1\right) B_{2(k+1)} }{(2(k+1))! }
x^{2k+1}. $$
Hence,  substituting it in (1.24), making the product of series and equating the coefficients in front of $t^{2n}$, we obtain
$$x\  p^\prime_n(x)=   \sum_{k=1}^{n}  {2n  \choose 2k-1  }  \frac{2^{2k}  -  1} {k }\ B_{2k } \  p_{n+1 -k} (x).\eqno(2.30)$$
In particular, dividing (2.30) by $x$ and passing $x$ to zero, we take into account the value $p^\prime_n(0)= -1,\ n \ge 1$ to derive the identity
$$\sum_{k=1}^{n}  {2n  \choose 2k-1  }  \frac{\left(2^{2k}  -  1\right)} {k }\ B_{2k }  = 1.\eqno(2.31)$$
Further,  employing  (1.21), we get from (2.30)
$$ \sum_{k=1}^{n} {2n +1 \choose 2k -1}
p_k(x) =    \sum_{r=1}^{n} {2n +1 \choose 2r } \sum_{k=1}^{r}  {2r  \choose 2k-1  }  \frac{\left(2^{2k}  -  1\right)} {k }\ B_{2k } \  p_{r+1 -k} (x).$$
Thus, multiplying  both sides of the latter equality by $e^{-x}$ and integrating over $\mathbb{R}_+$,  we use (2.20) to find the identity
$$ \sum_{r=1}^{n} \frac {E_{2r}}  {(2r-1)!  (2(n-r+1)) ! } =    2 \sum_{r=1}^{n}  \sum_{k=1}^{r}  \frac{\left(2^{2k}  -  1\right) B_{2k } \  E_{2(r-k +1)}} { (2k)! (2(n-r)+1)! (2(r-k)+1)!}.\eqno(2.32)$$

Returning to (2.17) and  letting $m=0$, we make a summation in  the right-hand side of (2.17) by $k$ from zero to $n$
$$\frac{n}{ 1- 2^{2n}} \sum_{k=0}^n {2n \choose 2k}  \int_0^\infty e^{-2x} p_{n-k}(x) p_{k}(x) {dx\over  x}$$
and employ  the binomial type identity (2.23)  to deduce
$$\sum_{k=0}^n {2n \choose 2k}  \int_0^\infty e^{-2x} p_{n-k}(x) p_{k}(x) {dx\over  x}=  \int_0^\infty e^{-x} p_{n}(x) {dx\over  x}.$$
Therefore,   we obtain the identity
$$ B_{2n}\sum_{k=0}^n {2n \choose 2k} = \frac{n \ }{ 1- 2^{2n}} \int_0^\infty e^{-x} p_n(x){dx\over  x}.\eqno(2.33)$$
An interesting question is to express the finite sum in the left-hand side of (2.33) in terms of the values $p_n^{\prime\prime} (0)$. In fact, differentiating two times in  (2.28) with respect to $x$, we let $x=0$ and use (2.12), (2.24) to  obtain
$$p^{\prime\prime}_n(0)= \sum_{k=2}^n  {2n  \choose 2k }q^{\prime\prime}_{k}(0)=  \sum_{k=0}^n  {2n  \choose 2k } - 2,$$
or,
$$  \sum_{k=0}^n  {2n  \choose 2k } = p^{\prime\prime}_n(0) +2,\ n \ge 1.\eqno(2.34)$$
Moreover,  differentiating two times in (2.30) and letting then $x=0$, we get the following recurrence relation for the values $p^{\prime\prime}_n(0)$
$$ p^{\prime\prime}_n(0)=  {1\over 2-n} \sum_{k=2}^{n-1}  {2n  \choose 2k-1  }  \frac{2^{2(n-k+1)}  -  1} {n-k+1 }\
B_{2(n-k+1) } \  p^{\prime\prime}_{k} (0),\ n\ge 1, \ n \neq 2, \eqno(2.35)$$
and $p^{\prime\prime}_2(0)=6.$

Nevertheless, we are able to calculate explicitly the left-hand side of (2.34) due to trigonometric and exponential series technique developed in \cite{gou} and where one can find a great collection of many such formulas.  Precisely, employing relation (3.7) in  Vol. 6, formula (3.7), it gives
$$  \sum_{k=0}^n  {2n  \choose 2k } =  2^{2n-1},  \  n \in \mathbb{N}\eqno(2.36)$$
and therefore,  $p^{\prime\prime}_n(0)=  2 \left(  2^{2(n-1)} -  1\right)$.   Consequently, identity (2.33) becomes
$$ B_{2n} = \frac{2n \ }{2^{2n}-  2^{4n} } \int_0^\infty e^{-x} p_n(x){dx\over  x}.\eqno(2.37)$$
Substituting the value of $p^{\prime\prime}_n(0)$ in (2.35), we obtain a possibly new identity
$$\sum_{k=2}^{n-1}  {2n  \choose 2k-1  }  \frac{\left(2^{2(n-k+1)}  -  1\right)\left( 2^{2(k-1)} -  1\right) } {n-k+1 }\
B_{2(n-k+1) } $$$$ =  \left(2^{2(n-1)} -  1\right) (2-n),\quad   n \ge 1.\eqno(2.38)$$
Meanwhile,
$$p_n(x)= \sum_{k=1}^n a_{n,k} x^k,\  n \ge 1,$$
where $a_{n,k}$ is defined by (1.15).  Thus, substituting it into (2.13) and (2.36), after calculation of the elementary Euler integral  and the use of (2.36), we find the following explicit formulas for the Bernoulli numbers, respectively,
$$B_{2n}= \frac{n}{ 1- 2^{2n}}  \sum_{k=1}^n   {1\over  2^k\   k}\sum_{r=0}^{k} {(-1)^r\over 2^r} \binom{k}{r}\sum_{j=0}^{k-r} {(-1)^j\over 2^j} \binom{k-r}{j} (r-j)^{2n} , \eqno(2.39)$$
$$B_{2n}= \frac{2n}{ 2^{2n} (1- 2^{2n}) }  \sum_{k=1}^n   {1\over k }\sum_{r=0}^{k} {(-1)^r\over 2^r} \binom{k}{r}\sum_{j=0}^{k-r} {(-1)^j\over 2^j} \binom{k-r}{j} (r-j)^{2n} \eqno(2.40)$$
and the equality of integrals
$$\int_0^\infty e^{-x} p_n(x){dx\over  x} =  2^{2n-1} \int_0^\infty e^{-2x} p_n(x){dx\over  x},\  n \ge 1.\eqno(2.41)$$
Concerning other identities and explicit formulas for Bernoulli's numbers see, for instance, a survey article \cite{gould} and in \cite{short}, \cite{ching}.   Further, substituting the right-hand side of (2.41) into (1.26),  it becomes
$$\zeta(2n)  =  \frac{ (-1)^n  \pi^{2n} }{2 (2^{2n}  - 1)(2n-1)!} \int_0^\infty e^{-x} p_n(x){dx\over x}, \ n \in \mathbb{N}.\eqno(2.42)$$
An explicit formula for the Euler numbers can be deduced similarly to (2.39), (2.40) with the use of the integral  (2.20).
Hence,  due to (1.15), we obtain for all $n \in \mathbb{N}$
 $$E_{2n}=  \sum_{k=1}^n  k! \sum_{r=0}^{k} {(-1)^r\over 2^r\ r!} \sum_{j=0}^{k-r} {(-1)^j(r-j)^{2n} \over 2^j\ j! (k-r-j)! } .$$
Other explicit formulas for Euler numbers see, for instance, in \cite{ming}.  The latter equality can give another characteristic of the Euler numbers.  In fact, we have
$$E_{2n}=  {d^{2n}\over d z^{2n} } \sum_{k=1}^n
   \sum_{r=0}^{k} {(-1)^r\over 2^r} \binom{k}{r}\sum_{j=0}^{k-r} {(-1)^j\over 2^j} \binom{k-r}{j} e^{z(r-j)}\ \Bigg|_{z=0} $$$$
   =    {d^{2n}\over d z^{2n} }  \sum_{k=1}^n  (1-\cosh z)^k \Bigg|_{z=0}.$$
Therefore, we find the formula
$$E_{2n} =    {d^{2n}\over d z^{2n} }  \   \frac{(1-\cosh z)(1 - (1-\cosh z )^{n})}{\cosh z} \Bigg|_{z=0}, \  n \ge 1.$$
Analogously, coefficients  (1.15) of the polynomial sequence $p_n(x)$ take the form
$$a_{n,k} = {(-1)^k 2^k\over k!}  {d^{2n}\over d z^{2n} }  \    \sinh^{2k}\left({z\over 2}\right) \Bigg|_{z=0},\   k=1,2, \dots, n .$$

In the meantime,  equality  (2.41) is quite important to derive  connection formulas  for the Bernoulli and Euler numbers. In fact, the integral in the left-hand side of (2.40) is calculated in \cite{yakpol} and we have
$$\int_0^\infty e^{-x} p_n(x){dx\over  x} =   -   \sum_{k=0}^{n-1} {2n -1 \choose 2k} E_{2k}, \  n \ge 1.$$
Consequently,  combining with  (2.13) and (2.37), we established the connection formula between the Bernoulli and Euler numbers.

{\bf Theorem 1}.  {\it  The following identity holds valid}
$$ B_{2n}= \frac{2n }{  2^{2n} (2^{2n}-  1)} \sum_{k=0}^{n-1} {2n -1 \choose 2k} E_{2k} ,\quad n \in \mathbb{N}.$$

Calling (2.42), we get an immediate

{\bf Corollary 1}.  {\it  For all $n \in \mathbb{N}$ one has }
$$ \zeta(2n)  = \frac{(-1)^{n+1} \pi^{2n}  }{2 (2^{2n}  - 1)}
\sum_{k=0}^{n-1} \frac{ E_{2k}}{(2k)! (2(n-k)-1) !}  ,\quad n \in \mathbb{N}.$$

Calling again (2.22), we differentiate through two times and let $x=0$. Hence with the use of (2.12) and (2.24) we derive a curious recurrence relation for  the Euler numbers.   Indeed, it has

{\bf Theorem 2}.  {\it  The following identity holds }
$$ E_{2n}= 1-   \sum_{k=0}^{n-1}  2^{2(n-k) -1} { 2n \choose 2k} E_{2k} ,\quad n \in \mathbb{N}.$$

As an application, we  announce  at the end of this section an interesting result about the structure of the Bernoulli numbers $B_{2n}$ and the rational values $\zeta(2n)/\pi^{2n}$ (see (2.10), (2.42)), which is a immediate  consequence of the Von Staudt- Clausen theorem \cite{von} about the fractional part of Bernoulli numbers and Fermat's Little theorem.

Precisely, it has

{\bf Theorem 3}.  {\it  The Bernoulli numbers $B_{2n}$ and Riemann zeta-values $\zeta(2n)$ satisfy the following properties, respectively,}
$$ 2 (  2^{2n}-  1) B_{2n}   \in \mathbb{Z},\quad n \in \mathbb{N},\eqno(2.43)$$
$$ 2 \left(2^{2n}  - 1\right)   {\zeta(2n) (2n-1)! \over \pi^{2n} }    \in \mathbb{Z},\quad n \in \mathbb{N}.\eqno(2.44)$$

Meanwhile identity (2.13) leads to

{\bf Corollary 2}.  {\it  For all $n \in \mathbb{N}$  }
$$2n  \int_0^\infty e^{-2x} p_{n}(x) {dx\over  x}  \in \mathbb{Z}.$$

\section {Riemann's  zeta-values}

Our main goal here is to establish certain identities, integral and series representations for the Riemann zeta function of positive argument.  Concerning  zeta-values at integers, as we could see in the previous section,  the Euler formula (2.10) gives a direct relationship  of $\zeta(2n),\ n \in \mathbb{N}$ with the Bernoulli numbers.  However,  similar formula for the values of zeta function at odd integers is unknown and probably does not exist.   Our attempts to find a finite relation between $\zeta(2n+1)$ are still  unsuccessful. Nevertheless, we will derive several integral and series representations,  related to these numbers and general positive numbers greater than one, involving our Sheffer's sequences of polynomials.   Some rapidly convergent series for $\zeta(2n+1)$ see, for instance, in \cite{sri}.

In fact, returning to  (1.28) and substituting  the modified Bessel function  by its representation (1.10), we employ the definition of the improper integral, integration by parts, the absolute and uniform convergence and the Riemann- Lebesgue lemma to make the change of the order of integration and motivate the following equalities
$$(-1)^n  (2n)!  \  \left(2^{2n+1} - 1\right)\   {\zeta(2n+1)\over (2 \pi)^{2n} } =  \lim_{N\to \infty} \int_0^N \tau \int_0^\infty  \int_0^\infty  e^{- 2x\cosh^2( u/2)} p_n(x) \cos(\tau u) {du dx d\tau \over  x} $$
$$=    \lim_{N\to \infty} \int_0^\infty  \int_0^\infty  e^{- 2x\cosh^2( u/2)} \sinh u \ p_n(x)\frac{1-\cos(Nu)}{u} du dx  $$
$$=  \int_0^\infty  \int_0^\infty  e^{- 2x\cosh^2( u/2)} \  p_n(x)\frac{\sinh u}{u} du dx = \int_0^1 \int_0^\infty  K_t(x) e^{-x}p_n(x) dx dt. $$
Consequently, we derived the identity for all $n \in \mathbb{N}$
$$(-1)^n  (2n)!  \  \left(2^{2n+1} - 1\right)\   {\zeta(2n+1)\over (2 \pi)^{2n} }  =
 \int_0^1 \int_0^\infty  K_t(x) e^{-x}p_n(x) dx dt.\eqno(3.1) $$
In the meantime, integrals (1.28), (3.1) have relationships with integrals, involving the Bernoulli polynomials owing
to the following representations proved in \cite{yakpol}
$$B_{2n+1}\left({1- t\over 2}\right)= -  {2n+1\over 2^{2n+1}\pi}\sin\pi t \ \int_0^\infty
K_{t}(x)e^{-x}p_n(x)dx,\    |t| < 1,\eqno(3.2)$$
$$B_{2n+1}\left({1- i\tau\over 2}\right)=  {2n+1\over
2^{2n+1}\ \pi i}\ \sinh\pi\tau  \ \int_0^\infty
K_{i\tau}(x)e^{-x}p_n(x)dx, \tau \in \mathbb{R}.\eqno(3.3)$$
Hence, integrating in (3.2) with respect to $t \in (0,1) $  and taking into account (3.1) we derive the identity
$$(-1)^{n+1}  (2n+1)!  \  \left(2-  2^{-2n}\right)\   {\zeta(2n+1)\over (2\pi)^{2n+1} }  =
\int_0^1 B_{2n+1}\left({1- t\over 2}\right)\frac{dt} { \sin\pi t },\   n\ge 1.\eqno(3.4)$$
Moreover, using (2.6), (2.7), (2.8),   after integration by parts with elementary substitutions in (3.4) and elimination of  the integrated terms, we write it in the form
$$(-1)^{n+1}  (2n)!  \  \left(2-  2^{-2n}\right)\   {\zeta(2n+1)\over 2^{2n+1} \pi^{2n}}  =
\int_0^{1/2} B_{2n}\left(t\right)\log\left(\cot \pi t  \right) dt,\   n\ge 1.\eqno(3.5)$$
One can find a similar identity, for instance, in \cite{ito}. The integral in (3.4) can be reduced  via properties for the Bernoulli polynomials to certain integrals considered recently in \cite{fur}.  Furthermore,  appealing to the addition formula for the Bernoulli polynomials \cite{erd}, Vol.  I  
$$B_n(x+y)=    \sum_{k=0}^{n}   { n \choose k} B_{k}(x) y^{n-k}\eqno(3.6)$$
the integral (3.4) can be represented as a linear combination of the moment integrals of $1/\sin t$, which we denote by $I_n$ following \cite{fur}
$$I_n= \int_0^{\pi/2} {t^n\over \sin t} dt,\  n \in \mathbb{N}.$$  
Hence,
$$\int_0^1 B_{2n+1}\left({1- t\over 2}\right)\frac{dt} { \sin\pi t }=  \int_0^{1/2}  B_{2n+1}\left({1- t\over 2}\right)\frac{dt} { \sin\pi t }+ \int_0^{1/2}  B_{2n+1}\left({t\over 2}\right)\frac{dt} { \sin\pi t } $$
$$= {1\over \pi}  \sum_{m=0}^{n}   (2\pi)^{2(m-n)-1} { 2n +1\choose 2m} \left[B_{2m}- B_{2m} \left({1\over 2}\right)\right] 
I_{2(n-m)+1} -   {2n+1\over (2\pi)^{2n+1}}I_{2n}.$$
But
$$B_{2m} \left({1\over 2}\right) =  - \left(1- 2^{1-2m}\right) B_{2m}.$$ 
Consequently, taking this value into account, we substitute the previous sum into  (3.4) and after simplifications arrive at the Ramanujan- type identity (cf. \cite{fur}) for zeta -values at odd integers $(n \in \mathbb{N}) $
$$(-1)^{n+1}   (2n+1)! \left(1-  2^{-2n-1}\right)   \zeta(2n+1) +  \left(n+{1\over 2}\right) I_{2n}  $$$$
 =  \sum_{m=0}^{n-1}   { 2n+1  \choose 2m+1 }   B_{2(n-m) } \pi^{2(n-m)-1} \left(2^{2(n-m)}- 1\right) I_{2m+1}.\eqno(3.7)$$ 
In particular, letting $n=1, 2$ and using the well-known formula $I_1= 2G $, where $G$ is the Catalan constant, we get, respectively, 
$${7\over 2} \zeta(3) + I_2=  2\pi G,$$
$$ I_4- {93\over 2} \zeta(5) =  2\pi ( I_3-  \pi^2 G).$$

{\bf Remark 1}.  We note, that the latter identities can be also obtained from corresponding equalities in Example I in \cite{fur}.

In the same manner one can obtain a finite sum representation of the zeta-values at odd integers in terms of the moment integrals of  $1/\sin^2 t$.  In this case our starting point will be Ito's identity \cite{ito}
$$(-1)^{n}  (2n)!  \   {\zeta(2n+1)\over (2\pi)^{2n} }  =
\int_0^{1} B_{2n}\left(t\right)\log\left(\sin \pi t  \right) dt,\   n\ge 1.\eqno(3.8)$$
Indeed, with the use of (2.8)  we have,
$$\int_0^{1} B_{2n}\left(t\right)\log\left(\sin \pi t  \right) dt =  \int_0^{1/2} B_{2n}\left(t\right)\log\left(\sin \pi t  \right) dt$$$$+
\int_0^{1/2} B_{2n}\left(1-t\right)\log\left(\sin \pi t  \right) dt= 2 \int_0^{1/2} B_{2n}\left(t\right)\log\left(\sin \pi t  \right) dt.$$
Hence, appealing to (2.4) in the right-hand side of the latter  equality and integrating twice by parts in the obtained integral, we  find
$$\int_0^{1} B_{2n}\left(t\right)\log\left(\sin \pi t  \right) dt  =   2 \sum_{m=0}^{n}  { 2n \choose 2m}  B_{2m}  \int_0^{1/2}  t^{2(n-m)} \log\left(\sin \pi t  \right) dt$$$$ -  2n \int_0^{1/2}  t^{2n-1} \log\left(\sin \pi t  \right) dt ={\pi^{-2n}\ M_{2n+1} \over 2n+1}    -  \sum_{m=0}^{n}  { 2n \choose 2m}  \frac{ B_{2m} \ \pi^{-2(n-m)-1} \  M_{2(n+1-m)}}{(2(n-m)+1)( n-m+1) } ,$$
where 
$$M_n=   \int_0^{\pi/2}  {t^{n} \over \sin^2 t} dt,\quad n \ge 2.$$
Thus combining with (3.8), we derived the identity (compare with (3.7))
$$(-1)^{n+1} (2(n+1))! \    2^{-2n-1} \zeta(2n+1)  +  (n+1) M_{2n+1} $$$$  =  
\sum_{m=0}^{n}  { 2(n+1)  \choose 2(m+1) }  B_{2(n-m)} \  \pi^{2(n-m)-1} \  M_{2(m+1)},\   n \ge 1.\eqno(3.9)$$
Appealing to relations (2.5.4.7) in \cite{prud},  Vol. I,  we have the values 
$$M_2= \pi \log 2,   \quad\quad    M_4=  {\pi^3\over 2}\log 2-  {9\pi\over 4} \zeta(3).$$
Therefore, letting $n=1$ in (3.9),  we get, for instance,
$$ {21\over 8}  \zeta (3) +  M_3= {3\pi^2 \over 4} \log 2.$$
Further, returning to (3.3),  we   multiply its  both sides by $\tau/\sinh\pi\tau $ and integrate over $\mathbb{R}_+$.  Hence
$$\int_0^\infty \tau  B_{2n+1}\left({1- i\tau \over 2}\right)\frac{d\tau} { \sinh\pi \tau }
=  {2n+1\over  2^{2n+1}\ \pi i}\int_0^\infty  \int_0^\infty \tau K_{i\tau}(x)e^{-x}p_n(x)dxd\tau $$$$
= {2n+1\over  2^{2n+1}\ \pi i}\int_0^1 \int_0^\infty  K_{t}(x)e^{-x}p_n(x)\  x  dx dt
,\   n\in \mathbb{N}_0.\eqno(3.10)$$
On the other hand, calling relation (1.20), we deduce from (3.1) and (3.10)
$$\sum_{k=0}^{n} {2n +1 \choose 2k}   {2^{2k+1} \over 2k+1} \int_0^\infty i \tau  B_{2k+1}\left({1- i\tau \over 2}\right)\frac{d\tau} { \sinh\pi \tau }$$$$ = {1\over \pi }  \sum_{k=0}^{n} {2n +1 \choose 2k}
\int_0^1 \int_0^\infty  K_{t}(x)e^{-x}p_k(x)\  x  dx dt$$$$ =  -  {1\over \pi } \int_0^1 \int_0^\infty  K_{t}(x)e^{-x}p_{n+1}(x) dx dt =  (-1)^{n}   (2(n+1))!  \  \left(2- 2^{- 2(n+1)}\right)\   {\zeta(2n+3)\over \pi^{2n+3} } $$
or
$$ \sum_{k=0}^{n} \frac {2^{2k}} {(2k+1)! (2(n-k)+1)! } \int_0^\infty i \tau  B_{2k+1}\left({1- i\tau \over 2}\right)\frac{d\tau} { \sinh\pi \tau } $$$$=   (-1)^{n}  (n+1) \  \left(2- 2^{- 2(n+1)}\right)\   {\zeta(2n+3)\over\pi^{2n+3} } .\eqno(3.11)$$
Hence, recalling  (3.6),  the integral in (3.11) can be rewritten as follows
$$2^{2k} \int_0^\infty i \tau  B_{2k+1}\left({1- i\tau \over 2}\right)\frac{d\tau} { \sinh\pi \tau } =  
 \sum_{m=0}^{k}  (-1)^{k-m+1}   { 2k +1\choose 2m}  \left( 2^{2m-1} - 1\right) B_{2m} $$$$\times  \int_0^\infty  \frac{\tau^{2(k-m+1)}\  d\tau} { \sinh\pi \tau }=   (2k+1)!  \sum_{m=0}^{k}  (-1)^{m+1} (m+1) \  \frac{\left( 2^{2(k-m)} - 2\right) \left( 2^{2(m +2)} - 2\right)B_{2(k-m)}} {(2(k-m))! } $$$$\times  { \zeta(2m+3)\over 
 (2\pi)^{2m+3} }.$$ 
Substituting the right-hand side of the latter equality in (3.11), we find the identity
$$ (-1)^{n}  (n+1) \  \left(2- 2^{- 2(n+1)}\right)\   {\zeta(2n+3)\over\pi^{2n+3} } =
 \sum_{k=0}^{n} \sum_{m=0}^{k}  \frac { (-1)^{m+1} (m+1)} {(2(n-k)+1)! } $$$$ \times   \frac{\left( 2^{2(k-m)} - 2\right) \left( 2- 2^{-2(m +1)}\right)B_{2(k-m)}} {(2(k-m))!  \pi^{2m+3} }  \zeta(2m+3)
 = \sum_{m =0}^{n}   (-1)^{m+1} (m+1) \left( 2- 2^{-2(m +1)}\right) $$$$\times  \frac{\zeta(2m+3)} { \pi^{2m+3}} 
\left(  \sum_{k=0}^{n-m}   \frac{\left( 2^{2k} - 2\right) B_{2k}} {(2k)! (2(n-m- k)+1)!   }\right) . \eqno(3.12)$$ 
It would be a great achievement to have here a finite recurrence relation for zeta- values at odd integers. However, unfortunately, this is not the case.  In fact, (3.12) yields for all $n \ge 1$
$$ \sum_{m =0}^{n-1}   (-1)^{m+1} (m+1) \left( 2- 2^{-2(m +1)}\right) \frac{\zeta(2m+3)} { \pi^{2m+3}} 
\left(  \sum_{k=0}^{n-m}   \frac{\left( 2^{2k} - 2\right) B_{2k}} {(2k)! (2(n-m- k)+1)!   }\right)= 0. \eqno(3.13)$$ 

{\bf Theorem 4}. {\it For all $n \in \mathbb{N}$ the following identity holds for Bernoulli numbers}
$$ \sum_{k=0}^{n}   { 2n+1  \choose 2k } \left( 2^{2k-1} - 1\right) B_{2k}   = 0.\eqno(3.14)$$

\begin{proof}  In fact, recalling (2.3), we see that (3.14) is equivalent to the equality
$$\sum_{k=0}^{n}   { 2n+1  \choose 2k }  2^{2k} B_{2k}   = 2n+1,$$
which yields  
$$\sum_{k=0}^{2n+1}   { 2n+1  \choose k }  2^{k-2n-1} B_{k}   = 0.$$
But this is true, because the left-hand side is equal (see (2.4), (2.8)) to $B_{2n+1} (1/2) = 0.$
\end{proof}
Theorem 4 says that all coefficients in front of zeta-values $\zeta(2m+3)$ in (3.14) are equal to zero.  Hence such kind of equalities  can be a source to obtain possibly new identities for Bernoulli numbers.    

Finally, we will get an integral representation of zeta-values at  positive numbers, which is a direct consequence of the formulas (1.3), (1.25), (1.26).  Precisely, it has 

{\bf Theorem 5}. {\it Let $\alpha > 1,\  [\alpha]$ be its integer part and $\{\alpha\}$ be its fractional part.  Then the following identities take place, when $[\alpha]$ is even or odd, respectively, }
$$\frac{2^\alpha - 1}{(2\pi)^{\alpha-1}}\  \Gamma(\alpha) \zeta(\alpha) = (-1)^{[\alpha]/2} 
\int_0^\infty \int_0^\infty \tau^{ \{\alpha\}} K_{i\tau}(x)\ e^{-x} p_{[\alpha]/2} (x){dx\over x},$$
$$\frac{2^\alpha - 1}{(2\pi)^{\alpha-1}}\  \Gamma(\alpha) \zeta(\alpha) = (-1)^{([\alpha]-1) /2} 
\int_0^\infty \int_0^\infty \tau^{ \{\alpha\}-1} K_{i\tau}(x)\ e^{-x} p_{([\alpha]+1)/2} (x){dx\over x}.$$

\bigskip
\centerline{{\bf Acknowledgments}}
\bigskip
The present investigation was supported, in part,  by the "Centro de  Matem{\'a}tica" of the University of Porto.

\vspace{5mm}

\noindent Semyon  Yakubovich\\
Department of  Mathematics,\\
Faculty of Sciences,\\
University of Porto,\\
Campo Alegre st., 687\\
4169-007 Porto\\
Portugal\\
E-Mail: syakubov@fc.up.pt\\

\end{document}